\documentclass[10pt,twosided]{article}
\usepackage[tbtags]{amsmath}
\usepackage{epsfig,amstext,amssymb,amsthm,latexsym}
\allowdisplaybreaks[4]
\usepackage{color}

\usepackage{amssymb,color}

\definecolor{c20}{rgb}{0.,0.7,0.}
\definecolor{c30}{rgb}{0.,0.,1.}
\definecolor{c40}{rgb}{1,0.1,0.7}
\definecolor{c50}{rgb}{1,0,0}
\definecolor{c60}{rgb}{0,0.9,0.1}

\def\rE#1{\textcolor{c20}{#1}}
\def\Hr#1{\textcolor{c20}{#1}}
\def\rE#1{#1}
\def\Hr#1{#1}

\def\xH#1{\textcolor{c20}{#1}}
\def\xH#1{#1}
\def\xx#1{\textcolor{c60}{#1}}
\def\xx#1{#1}

\def\yy#1{\textcolor{c60}{#1}}
\def\yy#1{#1}

\def\dE#1{#1}

\def\cJ#1{\textcolor{c50}{#1}}
\def\cJ#1{#1}

\def\cJI#1{\textcolor{c50}{#1}}
\def\cJI#1{#1}
\def\cL#1{\textcolor{c50}{#1}}
\def\cL#1{#1}
\def\cc#1{\textcolor{c50}{#1}}
\def\cc#1{#1}

\usepackage{amssymb,color}

\newcommand{\ve}{\varepsilon}

\newcommand{\abs}[1]{\lvert #1 \rvert}

\newcommand{\E}[1]{\mathbb{E}\left(#1\right)}

\newcommand{\pk}[1]{\mathbb{P} \left( #1 \right) }

\newcommand{\R}{\!I\!\!R}

\newcommand{\N}{\!I\!\!N}
\newcommand{\inr}{\in \R}

\newcommand{\inn}{\in \N}

\newcommand{\limit}[1]{\lim_{#1 \to   \infty}}

\newcommand{\BQN}{\begin{eqnarray}}
\newcommand{\EQN}{\end{eqnarray}}
\newcommand{\BQNY}{\begin{eqnarray*}}
\newcommand{\EQNY}{\end{eqnarray*}}

\newcommand{\BS}{\begin{sat}}
\newcommand{\ES}{\end{sat}}
\newcommand{\BT}{\begin{theo}}
\newcommand{\ET}{\end{theo}}
\newcommand{\BK}{\begin{korr}}
\newcommand{\EK}{\end{korr}}

\newcommand{\BD}{\begin{de}}
\newcommand{\ED}{\end{de}}

\newcommand{\BIT}{\begin{itemize}}
\newcommand{\EIT}{\end{itemize}}
\newcommand{\BDI}{\begin{description}}
\newcommand{\EDI}{\end{description}}

\newcommand{\BRM}{\begin{remarks}}
\newcommand{\ERM}{\end{remarks}}

\newcommand{\BEL}{\begin{lem}}
\newcommand{\EEL}{\end{lem}}

\newtheorem{theo}{Theorem}[section]
\newtheorem{sat}[theo]{Proposition}
\newtheorem{de}[theo]{Definition}
\newtheorem{lem}[theo]{Lemma}

\newtheorem{korr}[theo]{Corollary}
\newtheorem{remark}[theo]{Remark}
\newtheorem{remarks}[theo]{Remarks}

\newcommand{\nelem}[1]{{Lemma \ref{#1}}}

\newcommand{\netheo}[1]{{Theorem \ref{#1}}}

\newcommand{\prooftheo}[1]{ \textsc{Proof of Theorem} \ref{#1} }

\newcommand{\prooflem}[1]{\textsc{Proof of Lemma} \ref{#1}}

\newcommand{\COM}[1]{}

\newcommand{\QED}{\hfill $\Box$}

\def\Ga{\gamma}
\def\vp{\varepsilon}
\def\Sa{\sigma}
\def\SZV{\Sa_{Z_\vp}}
\def\rw{\rightarrow}
\newcommand{\mq}[1]{\quad\!\!\!\! \mathrm{#1}\quad\!\!\!\!}
\def\IF{\infty}
\def\piter{\mathcal{P}}
\newcommand{\expon}[1]{\exp\left(#1\right)}
\def\tildel{\tilde{\delta}}
\def\tildet{\tilde{t}_0}
\def\tildeY{\tilde{Y}}

\date{}

\def\mlH{\mathcal{H}}
\def\nonr{\nonumber}
\def\oo{(1+o(1))}
\def\asu{\ \ \text{as}\ u\rw\IF}

\begin{document}

\centerline{\bf \Large  \rE{On the} Supremum of $\gamma$-reflected Processes }
\bigskip
\centerline{\bf \Large  with Fractional Brownian Motion as Input}

\bigskip
\centerline{ Enkelejd Hashorva\footnote{Faculty of Business and Economics, University of
Lausanne, UNIL-Dorigny 1015 Lausanne, Switzerland, enkelejd.hashorva@unil.ch}, Lanpeng Ji$^a$ and Vladimir I. Piterbarg\footnote{Lomonosov
Moscow State University, 1 Leninskie Gory, 119234 Moscow, Russia,
piter@mech.math.msu.su}}

%\vskip -0.61 cm
%\centerline{\bf University of Lausanne \& Moscow State University}

  \centerline{\today{}}

\bigskip
{\bf Abstract:} Let $\{X_H(t), t\ge 0\}$ be a fractional Brownian motion with Hurst \rE{index} $H\in(0,\cJI{1]}$ \Hr{and
define} \cJI{a} $\gamma$-reflected process $W_\Ga(t)=X_H(t)-ct-\Ga\inf_{s\in[0,t]}\left(X_H(s)-cs \right)$, $t\ge0$ \Hr{with
$c>0,\gamma \in \xH{[0,1]}$ two given constants}.
In this paper we establish the exact tail asymptotic behaviour of \xH{$M_\gamma(T)=
\sup_{t\in [0,T]} W_\gamma(t)$}  for any $T\in (0,\IF]$.
Furthermore, we derive the exact tail asymptotic behaviour of the supremum of certain non-homogeneous \cJ{mean-zero} Gaussian random fields.
%It is surprising that the tail asymptotic behaviour of $M_\gamma(T)$ and $M_0(T)$ is up to some positive constant determined by $\gamma$ and $H$. \\% andthe \cJ{Piterbarg constant}.\\ %, the same as that of $M_0(T)$.



%%%%%%%%%%%%%%%%%%%%%%%%%%%%%%%%%%%5555555
\section{Introduction}
Let $\{X_H(t), t\ge0\}$ be  a standard fractional Brownian motion (fBm) with Hurst \rE{index}  $H\in(0,1]$, i.e., $X_H$ is
\cJ{a  $H$-self-similar Gaussian process with stationary increments,} and covariance  \xH{function}
\BQNY
Cov(X_H(t),X_H(s))=\frac{1}{2}(t^{2H}+s^{2H}-\mid t-s\mid^{2H}),\ t,s\ge0.
\EQNY
For two given constants  $\xH{c>0},\gamma\in [0,1]$ define  a new process $\{W_\gamma(t),t\ge 0\}$ by
\BQN
W_\Ga(t)=X_H(t)-ct-\Ga\inf_{s\in[0,t]}\left(X_H(s)-cs \right), \  \ t\ge0. \label{PCW}
\EQN
Throughout \xH{this} \rE{paper} $\{W_\gamma(t),t\ge 0\}$ \Hr{is referred to as a} {\it $\gamma$-reflected process} with fBm as input
 \xH{since} it reflects at rate $\gamma$ when reaching its minimum.\\
In queuing theory $W_1$ \xH{is the so called} {\it workload process}  (or {\it queue length process}) see e.g., Harrison (1985), Zeevi and Glynn (2000), Whitt (2002) \cL{and Awad and Glynn (2009)}; \Hr{alternatively one can} refer to $W_\gamma$ as a  generalized workload process with fBm as input.
\Hr{In} risk theory $W_\Ga$ can \Hr{be interpreted} as \Hr{a}  {\it claim surplus process} \Hr{since} the surplus process of an insurance \Hr{portfolio} can be defined by
$$U_\gamma(t)=u+ct-X_H(t)- \Ga\sup_{s\in[0,t]}\left(cs- X_H(s) \right)=u-W_\Ga(t),\ \ \  t\ge0$$
\Hr{for any nonnegative  initial reserve $u$.} \xH{In the literature, \rE{see e.g.,}  Asmussen and Albrecher (2010) the process $\{U_\gamma(t), t\ge 0\}$ is referred to as the risk process with tax payments of a {\it loss-carry-forward} type.}

This \Hr{contribution} is concerned with the tail asymptotic behaviour of \rE{the supremum} $M_\gamma(T)= \sup_{t\in[0,T]}W_\Ga(t)$ for $T \in (0,\IF]$, i.e.,
 we shall \Hr{investigate} the \Hr{rate} of convergence to 0 of
 $$\psi_{\Ga,T}(u):= \pk{M_\gamma(T)> u}$$
  as $u\to \IF$. The exact tail asymptotic behaviour of $M_\gamma(T)$ is known only for $\gamma=0$.  The case $T=\IF$ is already dealt with in H\"usler and Piterbarg (1999), whereas the case \cc{$T\in(0,\IF)$} has been \Hr{investigated} in D\c{e}bicki and Rolski (2002) and  D\c{e}bicki and Sikora (2011), see \Hr{our} \netheo{HPDS} in Appendix.
 Note in passing that  \cJ{ $M_1(t)\rw\IF$  almost surely as $t\rw\IF$ (e.g., Duncan and Jin (2008)),}
 therefore \Hr{we shall assume below} that $\gamma \in (0,1)$ when $T=\IF$. \\

The principal result of this paper is \xH{Theorem 1.1 below}, which establishes a unique asymptotic relationship between $\psi_{\Ga,T}(u)$
and $\psi_{0,T}(u)$ as $u\to \IF$ for any $T\in (0,\IF]$. Surprisingly, the following \yy{positive} constant %(commonly referred to as the Piterbarg constant)
\BQNY
\mathcal{P}_\alpha^a:=\underset{\rE{S}\rw\IF}\lim\mathcal{P}_\alpha^a[0,\rE{S}]=\underset{\rE{S}\rw\IF}\lim\mathcal{P}_\alpha^a[-\rE{S},0],\quad \alpha\in (0,2],\ a>0, %\ 0\le S_1<\IF,
\EQNY
where
\BQN\label{PS}
\mathcal{P}_\alpha^a[-S_1,S_2]=\E{ \exp\biggl(\sup_{t\in[-S_1,S_2]}\Bigl(\sqrt{2}B_{\alpha}(t)-(1+a)\abs{t}^{\alpha}\Bigr)\biggr)}\in (0,\IF),\quad
 0\le S_1,S_2<\IF,
 \EQN
 \cJ{with $\{B_{\alpha}(t),t\in\R\}$  a  fBm defined on $\R$ with Hurst \rE{index} $\alpha/2\in(0,1]$, }
determines the ratio
$$R_{\gamma,T}(u):=\frac{\psi_{\Ga,T}(u)}{\psi_{0,T}(u)}$$
for all $u$ large. Specifically, we have:

\BT \label{thnovel} For any $H,\gamma \in (0,1)$ and any $T\in (0,\IF]$
\BQN\label{novel}
\limit{u} R_{\gamma,T}(u)=\mathcal{M}_{H,\Ga,T} ,
\EQN
where $\mathcal{M}_{H,\Ga,T}=\piter_{2H}^{\frac{1-\Ga}{\Ga}}$ %if either $H< 1/2$ or
if $T=\IF$, and %. When $H\ge 1/2$ and
for $T< \IF$
$$
\mathcal{M}_{H,\Ga,T}=\left\{
            \begin{array}{ll}
\piter_{2H}^{\frac{1-\Ga}{\Ga}}  , & \hbox{if } \dE{H} <1/2 ,\\
\cJ{\piter_{1}^{\frac{\cJI{2}-\Ga}{\Ga}}}, & \hbox{if } \dE{H} =1/2,\\
1&  \hbox{if } \dE{H} >1/2.
              \end{array}
            \right.
$$
\ET
\cc{\rE{The} exact values of $ \mathcal{P}_\alpha^a$ are known only for $\alpha=1$ or $2$, namely,}
\BQN
 \mathcal{P}_1^a=1+\frac{1}{a}\ \ \text{and}\ \ \mathcal{P}_2^a=\frac{1}{2}\left(1+\sqrt{1+\frac{1}{a}}\right)
 \label{eqpp}
\EQN
see e.g., Piterbarg (1996) or D\c{e}bicki and Mandjes (2003). \cc{For general $\alpha\in(0,2)$, \rE{bounds for}  $ \mathcal{P}_\alpha^a$ are \rE{derived in D\c{e}bicki and Tabi\'{s} (2013)}.}

\rE{The asymptotic relation described by} \eqref{novel} is of relevance \Hr{for} theoretical models in queuing theory and insurance mathematics.
Moreover, a strong merit of \eqref{novel} is that its proof for the case $H> 1/2$ and $T=\IF$ is closely related with the exact tail asymptotics of
the supremum of \xH{certain} non-homogeneous \cJ{Gaussian random \xH{fields}}, a result which has not been known in the literature so far. Given the importance of that result for
extremes of Gaussian random fields, in the next section we present first an asymptotic expansion \Hr{of the tail of supremum of} certain non-homogeneous Gaussian random fields. We proceed then with the asymptotic formulas of $\psi_{\gamma, T}(u)$ for both cases $T=\IF$ and \cc{$T\in(0,\IF)$}. All proofs are relegated to Section 3 \Hr{followed} by some technical results displayed in Appendix.

\section{\Hr{Main} Results}
We prefer to state first our new result on the tail asymptotic behaviour of the supremum of certain non-homogeneous Gaussian fields,
since it \xH{is}
of theoretical importance going beyond \xH{the \yy{scope} of queuing and risk theory}.  We need to introduce some more notation starting with the well-known Pickands constant $\mathcal{H}_{\alpha}$ defined by
\cJI{$$\mathcal{H}_{\alpha}=\lim_{T\rightarrow\infty} \frac{1}{T}\mathcal{H}_{\alpha}[0,T],\ \ \alpha\in (0,2],$$
with $$\mathcal{H}_{\alpha}[0,T]=\E{ \exp\biggl(\sup_{t\in[0,T]}\Bigl(\sqrt{2}B_{\alpha}(t)-t^{\alpha}\Bigr)\biggr)}\in(0,\IF),\quad
\ T\in(0,\IF).$$}
 It is known that $\mlH_1=1$ and $\mlH_2={1}/{\sqrt{\pi}}$, see  Pickands (1969), Berman (1992),
 Piterbarg (1996), D\c{e}bicki (2002),  \cJI{Mandjes (2007)},  D\c{e}bicki and Mandjes (2011), D\c{e}bicki and
 \xx{Dieker and Yakir (2012)} for various properties of Pickands constant and its generalizations.\\
Throughout this paper, $x_+=\max(0,x)$ for any $x\inr$, and $\Psi(u), u\in\R$ denotes the \xx{survival} function
of the standard normal \xH{distribution $N(0,1)$}. \cc{Furthermore, we introduce the \rE{following constant}
$$\widetilde{\mathcal{P}}_\alpha^a=\lim_{\rE{S}\rw\IF}\mathcal{P}_\alpha^a[-\rE{S},\rE{S}],\ \  \alpha\in (0,2],\ a>0,$$
where $\mathcal{P}_\alpha^a[-\rE{S},\rE{S}]$ is given as in \eqref{PS}.}

We state next our first result.
\BT\label{ThmGPiter}
Let $S,T$ be two positive constants, and \cJI{let} $\{X(s,t),(s,t)\in\lbrack0,S]\times\lbrack0,T]\}$ be a zero-mean \cJ{Gaussian random field}, with \cc{standard deviation function} $\sigma(\cdot,\cdot)$ and correlation function $r(\cdot,\cdot,\cdot,\cdot)$.
Assume that $\sigma(\cdot,\cdot)$ \cc{attains its
unique maximum \xx{on} $\lbrack0,S]\times\lbrack0,T]$ at   $(s_0,t_0)$}, \cc{and further}
\BQN\label{eq:var}
\sigma(s,t)=1-b_{1}\abs{s-s_0}^{\beta}(1+o(1))-b_{2}(t-t_0)^2(1+o(1))- b_3\abs{(t-t_0)(s-s_0)}(1+o(1))\ \ \text{as }(s,t)\rightarrow(s_0,t_0)
\EQN
for some positive constants $b_i, i=1,2,3$ and $1<\beta<2$. Suppose further that
\BQN\label{eq:corr}
r(s,s^{\prime},t,t^{\prime})=1-(a_1|s-s^{\prime}|^{\beta}+a_2|t-t^{\prime}|^\beta%
)(1+o(1)) \text{\ \ \ as }(s,t),(s^{\prime
},t^{\prime})\rightarrow({s_0},{t_0})
\EQN
for some positive constants $a_i, i=1,2.$
If there \cL{exist two positive constants $G, \mu$ with $\mu\in(0,2]$} such that
\BQN\label{eq:Inc}
\E{(X(s,t)-X(s',t'))^2}\le G(|s-s^{\prime}|^{\xx{\mu}} +|t-t^{\prime}|^{\xx{\mu}})
\EQN
for any $(s,t), (s',t')\in\lbrack0,S]\times\lbrack0,T]$, then  %as $u\to \IF$
\BQN\label{eq:main}
\pk{\sup_{(s,t)\in\lbrack0,S]\times\lbrack0,T]}X(s,t)>u}=\widehat{\mathcal{P}_{\beta}^{b_1/a_1}}\ \widehat{I}_2 \ \mathcal{H}_\beta \frac{\sqrt{\pi}a_2^{1/\beta}}{2\sqrt{b_2}} u^{2/\beta-1}\Psi(u)  (1+o(1))\asu, %\int_{T_2^\infty} e^{-b_2 t^2}dt
\EQN
\cc{where the constant $\widehat{\mathcal{P}_{\beta}^{b_1/a_1}}$ is equal to $\widetilde{\mathcal{P}}_{\beta}^{b_1/a_1}$ if $s_0\in(0,S)$ and equal to $\mathcal{P}_{\beta}^{b_1/a_1}$ if $s_0=0$ or $S$, and $\widehat{I}_2$ is equal to 2  if $t_0\in(0,T)$ and equal to 1 if $t_0=0$ or $T$.}
%with $T_1^{\infty}=\underset{{x\rw\infty}}\lim x([0,S]-s_0)$ and $T_2^{\infty}=\underset{{x\rw\infty}}\lim x([0,T]-t_0)$ (e.g., $\underset{x\rw\infty}\lim\ x[0,1]=[0,\IF)$).
\ET
\begin{remark}\label{b3n}
\cJI{The claim of \netheo{ThmGPiter} still holds for $b_3<0$ satisfying $b_2+b_3/2>0$ \rE{which can be shown by utilising \netheo{ThmGPiter} and
noting that}
\BQNY
1-\sigma(s,t)\ge\frac{b_{1}}{2}\abs{s-s_0}^{\beta}(1+o(1))+\left(b_{2}+\frac{b_3}{2}\right)(t-t_0)^2(1+o(1))%- b_3\abs{t-t_0}\abs{s-s_0}(1+o(1))
\EQNY
as $(s,t)\rightarrow(s_0,t_0)$.} \xx{Note  that \eqref{eq:main} does not depend on the value of the constant $b_3$.}
\end{remark}

Next, we return to our principal problem deriving below the exact asymptotic behaviour of $\psi_{\gamma, T}(u)$ as $u \to \IF$.
Although \cJ{the limit of} the ratio $R_{\gamma,T}(u)$ \cL{as $u\rw\IF$} remains constant, both cases \cc{$T= \IF$ and $T\in(0,\IF)$} are very different and will \xH{therefore}
be dealt with separately. We \Hr{shall analyse} first the case $T=\IF$. \\

\yy{Below, we} set $Y_u(s,t):=\frac{X_H(ut)-\gamma X_H(us)}{(1+ct-c\Ga s)u^H}$, and then  write
\BQN
\psi_{\Ga,\IF}(u)&=&\pk{\sup_{t\ge 0}\left(X_H(t)-ct-\gamma\inf_{s\in[0,t]}(X_H(s)-cs)  \right)>u}\nonumber\\
%&=&\pk{X_H(t)-\gamma X_H(s)+c\gamma s -ct>u, \mq{for}\mq{some}s,t\mq{such}\mq{that}0\le s \le t<\IF}\nonumber\\
&=&\pk{\sup_{0\le s \le t<\IF} Y_u(s,t)>u^{1-H}}.
\label{Infruin}
\EQN
The above alternative formula for $\psi_{\Ga, \IF}(u)$ together with \netheo{ThmGPiter} and \nelem{lemVarY} is crucial for the derivation \cJ{of} the tail asymptotic behaviour of $M_\gamma(\IF)$. % given in \netheo{th2} below.

\BT\label{th2}  We have, for $H, \Ga \in(0,1)$
\BQN\label{t2}
\psi_{\Ga,\IF}(u)=\mathcal{W}_{H}(u)\Psi\left(\frac{c^H u^{1-H}}{H^H (1-H)^{1-H}}\right)(1+o(1))\ \ \text{as}\  u\rw\infty,
\EQN
where
$$
\mathcal{W}_{H}(u)=
 2^{\frac{1}{2}-\frac{1}{2H}}\frac{\sqrt{\pi}}{\sqrt{H(1-H)}} \mathcal{H}_{{2H}}\piter_{2H}^{\frac{1-\Ga}{\Ga}}\left(\frac{c^H u^{1-H}}{H^H (1-H)^{1-H}}\right)^{1/H-1}.
$$

\ET
%%%%%%%%%%%%%%%555555
\def\NN{\mathcal{N}}
%%%%%%%%%%55555
\begin{remark}
 If $H=1$, then  $X_H(t)=\NN t$ with $\NN$ \yy{a} standard normal random variable (i.e., $\NN$ \yy{is $N(0,1)$ distributed}).
 Consequently, for any $c>0$ \rE{and $\gamma \in (0,1)$}
\BQNY
\psi_{\Ga,\IF}(u)&=&\pk{\sup_{0\le s\le t<\IF}\frac{(t-\Ga s)\NN }{1+c(t-\Ga s)}>1}
= \Psi(c) =\pk{\sup_{t\ge 0}\Bigl(\NN t-ct\Bigr)>u}=\psi_{0,\IF}(u)
\EQNY
holds for all $u\ge0.$
%ii) In view of \eqref{eq:HP} in \netheo{HPDS} in Appendix the claim of \netheo{thnovel} for $T=\IF$ follows immediately by \eqref{t2}.
\end{remark}

{\bf Example}: Consider the case of the $\gamma$-reflected process with Brownian motion as input, i.e., $H=1/2$. It is well-known that
\BQNY
\psi_{0,\IF}(u)=\pk{\sup_{t\in[0,\IF)}\Bigl( B_{1}(t)-ct \Bigr)>u}=e^{-2cu}, \ \ u\ge0.
\EQNY
Further, \xH{for this case} \netheo{th2} \cL{together with} \eqref{eqpp} \rE{imply}
$$\psi_{\Ga,\IF}(u)=\frac{2\sqrt{2\pi}c^{1/2}}{1-\Ga}u^{1/2}\Psi\left(2c^{1/2} u^{1/2}\right)(1+o(1))$$
as $u\rw\IF.$
Therefore
$$
\psi_{\Ga,\IF}(u)={\frac{1}{1-\Ga}}\psi_{0,\IF}(u)(1+o(1)) \ \ \text{as}\ u\rw\IF,
$$
which \cc{also follows from the following} \rE{identity (see e.g., Asmussen and Albrecher (2010), Albrecher et al.\ (2013))}
$$
\psi_{\Ga,\IF}(u)+(1-\psi_{0,\IF}(u))^{\frac{1}{1-\Ga}}=\xH{1} , \ \ \forall u\ge0.
$$
%derived in Albrecher et al. (2008).} %is the asymptotic tax identity given in \eqref{eq:AH}.

We conclude this section with an explicit asymptotic expansion for $\psi_{\gamma, T}(u)$ with $T \in (0,\IF)$.
For any $u\ge 0$
\BQNY
\psi_{\Ga,T}(u)&=&\pk{\sup_{t\in[0,T]}\left(X_H(t)-ct-\gamma\inf_{s\in[0,t]}(X_H(s)-cs)  \right)>u}\nonumber\\
%&=&\pk{X_H(t)-\gamma X_H(s)+c\gamma s -ct>u, \mq{for}\mq{some}s,t\mq{such}\mq{that}0\le s \le t\le T}\nonumber\\
&=&\pk{\sup_{0\le s \le t\le T}(Z(s,t)-(ct -c\gamma s))>u},
\EQNY
where $Z(s,t):=X_H(t)-\gamma X_H(s)$. It follows that the variance function of $Z(s,t)$ is given by
\BQNY
\xx{V_Z}^2(s,t)=\E{(X_H(t)-\gamma X_H(s))^2}=(1-\gamma)t^{2H}+(\Ga^2-\Ga)s^{2H}+\Ga(t-s)^{2H}.
\EQNY
Clearly, $\xx{V_Z}(s,t)$ attains its unique  maximum \cc{\yy{on} the set $A:=\{(s,t):0\le s \le t\le T\}$ }  at $(s_0,t_0)=(0,T)$. This fact is crucial for our last result stated below.

\BT\label{th1}  For any \cJI{$T\in(0,\IF)$} and $\cJI{H}, \Ga\in(0,1],$ we have
\BQN
\psi_{\Ga,T}(u)=\mathcal{D}_{H,\Ga} \left(\frac{u+cT}{T^H}\right)^{\frac{(1-2H)_+}{H}}\Psi\left(\frac{u+cT}{T^H}\right)(1+o(1)) \ \ \text{as}\ u\rw\IF,\label{Zvp}
\EQN
where
$$
\mathcal{D}_{H,\Ga}=\left\{
            \begin{array}{ll}
2^{-\frac{1}{2H}}H^{-1}\mathcal{H}_{{2H}}\piter_{2H}^{\frac{1-\Ga}{\Ga}}  , & \hbox{if } \dE{H} <1/2 ,\\
\frac{4}{2-\Ga}, & \hbox{if } \dE{H} =1/2,\\
1&  \hbox{if } \dE{H} >1/2,
              \end{array}
            \right.\ \Ga\in(0,1),\ \
\mathcal{D}_{H,1}=\left\{
            \begin{array}{ll}
2^{-\frac{1}{H}}H^{-2}(\mathcal{H}_{{2H}})^2 , & \hbox{if } \dE{H} <1/2 ,\\
4, & \hbox{if } \dE{H} =1/2,\\
1&  \hbox{if } \dE{H} >1/2.
              \end{array}
            \right.
$$
\ET
\COM{The above In v%Note that the finite-time ruin probability without tax for fBm model was given in Debicki and Sikora (2011).
Similar to the infinite-time case, in the following corollary we give the tax-identity for the finite-time case. %ruin probability with tax and the finite-time ruin probability without tax.
\BK
We obtain, for any $H\in(0,1]$
\BQN\label{eq:FAH}
\psi_{\Ga,T}(u)=\mathcal{M}_{H,\Ga} \psi_{0,T}(u)(1+o(1)) \ \ \text{as}\ u\rw\IF,
\EQN
where
$$
\mathcal{M}_{H,\Ga}=\left\{
            \begin{array}{ll}
\piter_{2H}^{\frac{1-\Ga}{\Ga}}  , & \hbox{if } \dE{H} <1/2 ,\\
\frac{2}{2+\Ga}, & \hbox{if } \dE{H} =1/2,\\
1&  \hbox{if } \dE{H} >1/2.
              \end{array}
            \right.
$$
\EK
}

\section{Proofs }
%%%%%%%%%555555
\def\CC{\mathbb{C}}
%%%%%%%5555555
In this section, we give proofs of all the results. Hereafter the positive constant $\CC$ may
be different from line to line. Furthermore, a mean-zero \Hr{Gaussian} process (\cJI{or \Hr{a} random \Hr{Gaussian} field}) $\{\xi(t),t\ge0\}$ with a bar denotes the corresponding standardized process (\cJI{or random field}), i.e., $\overline{\xi}(t)=\xi(t)/\sqrt{\E{\xi(t)}^2}$.

\prooftheo{thnovel}  \cJI{In view of \eqref{eq:HP} in \netheo{HPDS} in Appendix the claim  for the case $T=\IF$ follows immediately from \netheo{th2}. Further, \Hr{by} combining the result of \netheo{th1} with that of \eqref{eq:DS} in \netheo{HPDS} \Hr{we} establishe the claim for the case \cc{$T\in(0,\IF)$}.} \QED

\bigskip
\prooftheo{ThmGPiter}
%%%%%%%%%5555
\def\delu{\delta(u)}
%%%%%%%%%55555
%In this Appendix we give a rigorous proof \netheo{ThmGPiter}.
Set $\eta(s,t)=X(s+s_0,t+t_0), (s,t)\in [-s_0,S-s_0]\times[-t_0,T-t_0]$. It follows that the standard deviation  function $\sigma_\eta(s,t)$ of $\eta(s,t)$  attains its unique maximum   \rE{equal to} 1 \cc{on $[-s_0,S-s_0]\times[-t_0,T-t_0]$}  at $(0,0)$. Further \eqref{eq:var} and \eqref{eq:corr} are valid for the standard deviation  function $\sigma_\eta$ and the correlation function $r_\eta$ with $(s_0,t_0)$ replaced by $(0,0)$. Moreover, \eqref{eq:Inc} is established for the random field $\eta$ over $[-s_0,S-s_0]\times[-t_0,T-t_0]$. There are \Hr{nine} cases to be considered depending on whether \xx{$0$ is} \yy{an} inner \cL{point or \cc{a} boundary point} of $[-s_0,S-s_0]$ \xx{or} $[-t_0,T-t_0]$. We \Hr{investigate} next on the case that $(s_0,t_0)=(0,0)$, and thus $\eta=X$. The other cases can be  \Hr{analysed} with the same argumentations.

\def\Du{\xH{\tilde \Delta_u}}
{\Hr{In the light of} Theorem 8.1 in Piterbarg (1996) (or Theorem 8.1 in Piterbarg (2001)) for $u$ sufficiently large (set
$\delta(u)=\ln u/u, \Du=[0,\delu]\times[0,\delu]$)}
\BQN\label{eq:neg}
\pk{\sup_{(s,t) \in [0, S]\times[0,T]/\Du } X(s,t)>u}\le \CC \ \cL{u^{4/\mu}} \exp\left(-\frac{u^2}{2-C\delu^2}\right)
\EQN
holds for some positive constant \Hr{$\CC$} not depending on $u$. Next we \yy{analyse}
\BQNY
\pk{\sup_{(s,t)\in\Du}X(s,t)>u}
\EQNY
\def\xiS{ \xH{\tilde \xi(s,t)}}
\Hr{as $u\to \IF$}, which has the same asymptotic \xx{behaviour} as (set $\xiS=\frac{\xi(s,t)}{(1+b_1s^\beta)(1+b_2t^2+b_3ts)}$)
\BQNY
\pi(u)=\pk{\sup_{(s,t)\in\Du} \xiS>u} \quad \text{as  } u\to \IF,
\EQNY
where $\{\xi(s,t), s,t\ge0\}$
is a \Hr{mean-zero} \cJ{Gaussian random field} with covariance function given by
$$
r_\xi(s,t)=\exp(-a_1s^\beta-a_2t^\beta), \quad s,t\ge0.
$$
\Hr{For simplicity, we shall} assume that $a_1=a_2=1.$ The general case can be \Hr{analysed} by \cc{rescaling} \Hr{the time}.
It follows from Lemma 6.1 \yy{in} Piterbarg (1996) that
\BQN\label{eq:piter6.1}
\pk{\underset{[0, u^{-2/\beta}S]\times[0,u^{-2/\beta}T]}{\sup_{(s,t)\in}}\frac{\xi(s,t)}{1+b_1s^\beta}>u}=\mathcal{P}_{\beta}^{b_1}[0,S]\mathcal{H}_\beta[0,T]\Psi(u)(1+o(1))\ \ \text{as}\  u\rw\IF.
\EQN
Since $\beta \yy{\in (1,2)}$, for any positive constant $S_1$, we can divide the interval $[0,\delu]$ into several sub-intervals of length $S_1u^{-2/\beta}$.  Specifically, let for $S_1,S_2>0$
%%%%%%%%%55
\def\det{\triangle}
\def\wdet{\widetilde{\triangle}}
\def\H{\mathcal{H}}
%%%%%%%
\BQNY
\det_0^i=u^{-2/\beta}[0,S_i], \ \ \ \det_k^i=u^{-2/\beta}[k S_i,(k +1)S_i], \ \ k \inn, \ i=1,2
\EQNY
\Hr{and let} further
$$
h_i(u)=\lfloor S_i^{-1}u^{\frac{2}{\beta}-1} \ln u \rfloor+1,\ \ i=1,2,\ u>0.
$$
Bonferroni inequality \Hr{yields}
\BQNY
\pi(u)&\le&
%\sum_{k_1=0}^{h_1(u)}\sum_{k_2=0}^{h_2(u)}\pk{\underset{(s,t)\in\det_{k_1}^1\times\det_{k_2}^2}\sup
%\xiS%\frac{\xi(s,t)}{(1+b_1s^\beta)(1+b_2t^2+b_3ts)}
%>u}\\
%&=&
\sum_{k_2=0}^{h_2(u)}\pk{\underset{(s,t)\in\det_{0}^1\times\det_{k_2}^2}\sup
\xiS %\frac{\xi(s,t)}{(1+b_1s^\beta)(1+b_2t^2+b_3ts)}
>u}
%\\
%&&
+\sum_{k_1=1}^{h_1(u)}\sum_{k_2=0}^{h_2(u)}\pk{\underset{(s,t)\in\det_{k_1}^1\times\det_{k_2}^2}\sup
\xiS%\frac{\xi(s,t)}{(1+b_1s^\beta)(1+b_2t^2+b_3ts)}
>u}\\
&=:&I_1(u)+I_2(u)
\EQNY
and
\BQNY
\pi(u)&\ge& \sum_{k_2=0}^{h_2(u)-1}\pk{\underset{(s,t)\in\det_{0}^1\times\det_{k_2}^2}\sup
\xiS%\frac{\xi(s,t)}{(1+b_1s^\beta)(1+b_2t^2+b_3ts)}
>u}\\
%&=&\sum_{k_2=0}^{h_2(u)}\pk{\underset{(s,t)\in\det_{0}^1\times\det_{k_2}^2}\sup\frac{\xi(s,t)}{(1+b_2s^\beta)(1+b_1t^2+b_3ts)}>u}\\
&&- \sum_{0\le i<j\le h_2(u)-1}
\pk{\underset{(s,t)\in\det_{0}^1\times\det_{i}^2}\sup\xiS%\frac{\xi(s,t)}{(1+b_1s^\beta)(1+b_2t^2+b_3ts)}
>u,\underset{(s,t)\in\det_{0}^1\times\det_{j}^2}\sup\xiS %\frac{\xi(s,t)}{(1+b_1s^\beta)(1+b_2t^2+b_3ts)}
>u}\\
&=:&J_1(u)- J_2(u).
\EQNY
Next we \rE{calculate} the \cJI{required} asymptotic bounds \rE{for}
$I_1(u)$ and $J_1(u)$ and show that
\BQN\label{eq:double0}
I_2(u)=J_2(u)(1+o(1))=o(I_1(u))\ \ \ \text{as}\  u\rw\IF,\ S_i\rw\IF, i=1,2.
\EQN
We derive that
\BQNY
J_1(u)&=& \sum_{k_2=0}^{h_2(u)-1}\pk{\underset{(s,t)\in\det_{0}^1\times\det_{k_2}^2}\sup\xiS%\frac{\xi(s,t)}{(1+b_1s^\beta)(1+b_2t^2+b_3ts)}
>u}\nonumber\\
&\ge&\sum_{k_2=0}^{h_2(u)-1}\pk{\underset{(s,t)\in\det_{0}^1\times\det_{k_2}^2}\sup\frac{\xi(s,t)}{1+b_1s^\beta}
>u(1+b_2((k_2+1)S_2u^{-2/\beta})^2+b_3((k_2+1)S_2u^{-2/\beta})(S_1u^{-2/\beta}))}.
\EQNY
In view of \eqref{eq:piter6.1}
\BQN\label{eq:J_1}
J_1(u)&\ge& \mathcal{P}_{\beta}^{b_1}[0,S_1]\mathcal{H}_\beta[0,S_2]\frac{1}{\sqrt{2\pi}}\sum_{k_2=0}^{h_2(u)-1}
\frac{1}{u(1+b_2((k_2+1)S_2u^{-2/\beta})^2+b_3((k_2+1)S_2u^{-2/\beta})(S_1u^{-2/\beta}))}\nonumber \\ &&\exp\left(-\frac{u^2(1+b_2((k_2+1)S_2u^{-2/\beta})^2+b_3((k_2+1)S_2u^{-2/\beta})(S_1u^{-2/\beta}))^2}{2}\right)\oo\nonumber\\
%&=& \mathcal{P}_{\beta}^{b_1}[0,S_1]\mathcal{H}_\beta[0,S_2] \Psi(u)      %\frac{1}{\sqrt{2\pi}u}\exp\left(-\frac{u^2}{2}\right)
%\sum_{k_2=0}^{h_2(u)-1}\exp\left(-b_2((k_2+1)S_2u^{1-2/\beta})^2-b_3 S_1 S_2((k_2+1)u^{2-4/\beta})\right)\oo\nonumber\\
%&=& \mathcal{P}_{\beta}^{b_1}[0,S_1]\frac{\mathcal{H}_\beta[0,S_2]}{S_2} \Psi(u)   u^{2/\beta-1}   %\frac{1}{\sqrt{2\pi}u}\exp\left(-\frac{u^2}{2}\right)
%\sum_{k_2=0}^{h_2(u)-1}\exp\left(-b_1((k_2+1)S_2u^{1-2/\beta})^2-b_3 S_1 S_2((k_2+1)u^{2-4/\beta})\right)S_2u^{1-2/\beta}\nonumber\\
&=&\mathcal{P}_{\beta}^{b_1}[0,S_1]\frac{\mathcal{H}_\beta[0,S_2]}{S_2} \Psi(u)   u^{2/\beta-1} \int_0^\IF e^{-b_2x^2}dx\oo
\EQN
as $u\rw\IF,$
where in the last equation we used the facts that, as $u\rw\IF$
$$
h_2(u)\rw\IF,\ \ h_2(u)u^{1-2/\beta}\rw\IF,\ \ h_2(u)u^{2-4/\beta}\rw 0.
$$
Similarly
\BQN\label{eq:I_1}
I_1(u)
&\le&
%\sum_{k_2=0}^{h_2(u)-1}\pk{\underset{(s,t)\in\det_{0}^1\times\det_{k_2}^2}\sup\frac{\xi(s,t)}{1+b_1s^\beta}
%>u(1+b_2(k_2S_2u^{-2/\beta})^2)}\nonumber\\
%&=&
\mathcal{P}_{\beta}^{b_1}[0,S_1]\frac{\mathcal{H}_\beta[0,S_2]}{S_2} \Psi(u)   u^{2/\beta-1} \int_0^\IF e^{-b_2x^2}dx\oo
\EQN as $u\rw\IF$.
Moreover, \eqref{eq:double0} can be shown as in Piterbarg (1996). Specifically
\BQNY
I_2(u)&=&\sum_{k_1=1}^{h_1(u)}\sum_{k_2=0}^{h_2(u)}\pk{\underset{(s,t)\in\det_{k_1}^1\times\det_{k_2}^2}\sup
\xiS%\frac{\xi(s,t)}{(1+b_1s^\beta)(1+b_2t^2+b_3ts)}
>u}\\
&\le&\sum_{k_1=1}^{h_1(u)}\sum_{k_2=0}^{h_2(u)}\pk{\underset{(s,t)\in\det_{k_1}^1\times\det_{k_2}^2}\sup \xi(s,t) >u(1+b_1(k_1S_1u^{-2/\beta})^\beta+b_2(k_2S_2u^{-2/\beta})^2)}. %+b_3(k_1S_1u^{2/\beta})(k_2S_2u^{2/\beta})
\EQNY
Similar argumentations as in \eqref{eq:J_1} yield
\BQNY
%&\overset{\eqref{eq:piter6.1}}\sim&\H_\beta[0,S_1]\H_\beta[0,S_2]\Psi(u)(S_2^{-1}u^{2/\beta-1})\sum_{k_1=1}^{h_1(u)}\sum_{k_2=0}^{h_2(u)}
%\exp\left(-b_1(k_2S_2u^{1-2/\beta})^2-b_2 (k_1S_1)^\beta\right)S_2u^{1- 2/\beta},\\
I_2(u)&\le&\H_\beta[0,S_1]\H_\beta[0,S_2]\Psi(u)(S_2^{-1}u^{2/\beta-1})\int_0^\IF e^{-b_2x^2}dx\sum_{k_1=1}^{h_1(u)}
\exp\left(-b_1 (k_1S_1)^\beta\right)\oo
\EQNY as $u\rw\IF$.
Further
\BQNY
J_2(u)&=&\sum_{0\le i<j\le h_2(u)-1}\pk{\underset{(s,t)\in\det_{0}^1\times\det_{i}^2}\sup\xiS%\frac{\xi(s,t)}{(1+b_1s^\beta)(1+b_2t^2+b_3ts)}
>u,
\underset{(s,t)\in\det_{0}^1\times\det_{j}^2}\sup\xiS%\frac{\xi(s,t)}{(1+b_1s^\beta)(1+b_2t^2+b_3ts)}
>u}=:\Sigma_1(u)+\Sigma_2(u),
\EQNY
where $\Sigma_1(u)$ is the sum \rE{over indices}  $j=i+1$, and similarly $\Sigma_2(u)$ is the sum over indices  $j>i+1$.
Let
$$
B(i,S_2,u)=u(1+b_2(iS_2u^{-2/\beta})^2), \ \ i\inn,\ S_2>0,\ u>0.
$$
It follows that
\BQNY
\Sigma_1(u)&\le&\sum_{0\le i \le h_2(u)-1}\pk{\underset{(s,t)\in\det_{0}^1\times\det_{i}^2}\sup\frac{\xi(s,t)}{1+b_1s^\beta}>B(i,S_2,u),
\underset{(s,t)\in\det_{0}^1\times\det_{i+1}^2}\sup\frac{\xi(s,t)}{1+b_1s^\beta}>B(i,S_2,u)}
\EQNY
and
\BQNY
&&\pk{\underset{(s,t)\in\det_{0}^1\times\det_{i}^2}\sup\frac{\xi(s,t)}{1+b_1s^\beta}>B(i,S_2,u),
\underset{(s,t)\in\det_{0}^1\times\det_{i+1}^2}\sup\frac{\xi(s,t)}{1+b_1s^\beta}>B(i,S_2,u)}\\
&&=\pk{\underset{(s,t)\in\det_{0}^1\times\det_{0}^2}\sup\frac{\xi(s,t)}{1+b_1s^\beta}>B(i,S_2,u)}+\pk{
\underset{(s,t)\in\det_{0}^1\times\det_{1}^2}\sup\frac{\xi(s,t)}{1+b_1s^\beta}>B(i,S_2,u)}\\
&&-\pk{\underset{(s,t)\in\det_{0}^1\times(\det_{0}^2\cup\det_1^2)}\sup\frac{\xi(s,t)}{1+b_1s^\beta}>B(i,S_2,u)}.
\EQNY
Therefore, using the same reasoning as \eqref{eq:J_1}, we conclude that % \eqref{eq:sig1}
\BQN\label{eq:sig1}
\limsup_{u\rw\IF}\frac{\Sigma_1(u)}{\Psi(u)   u^{2/\beta-1}}&\le&
\mathcal{P}_{\beta}^{b_1}[0,S_1] \frac{2\mathcal{H}_\beta[0,S_2]-\H_\beta[0,2S_2]}{S_2}\int_0^\IF e^{-b_2x^2}dx .
\EQN
Further
\BQNY
\Sigma_2(u)
%&\le&\sum_{i=0}^{ h_2(u)-1}\sum_{j\ge2}\pk{\underset{(s,t)\in\det_{0}^1\times\det_{0}^2}\sup \xi(s,t) >B(i,S_2,u),
%\underset{(s,t)\in\det_{0}^1\times\det_{j}^2}\sup \xi(s,t) >B(i,S_2,u)}\\
&\le&\sum_{i=0}^{ h_2(u)-1}\sum_{j\ge2}\pk{\underset{(s',t')\in\det_{0}^1\times\det_{j}^2}{\sup_{(s,t)\in\det_{0}^1\times\det_{0}^2}} \zeta(s,t,s',t') >2B(i,S_2,u)},
\EQNY
where
$$
\zeta(s,t,s',t')=\xi(s,t)+\xi(s',t'),\ \ s,s',t,t'\ge0.
$$
\Hr{Now}, for $u$ sufficiently large
$$
2\le\E{(\zeta(s,t,s',t'))^2}=4-2(1-r(\abs{s-s'},\abs{t-t'}))\le 4- ((j-1)S_2)^\beta u^{-2}
$$
for any $ (s,t)\in\det_{0}^1\times\det_{0}^2, (s',t')\in\det_{0}^1\times\det_{j}^2$.
Thus, using similar argumentations as in Lemma 6.3 of Piterbarg (1996), we conclude that
\BQN\label{eq:sig2}
\limsup_{u\rw\IF}\frac{\Sigma_2(u)}{\Psi(u)   u^{2/\beta-1}}&\le& \CC\
(\mathcal{H}_{\beta}[0,S_1])^2  S_2 \sum_{j\ge1}\exp\left(-\frac{1}{8}(jS_2)^\beta\right).
\EQN

Consequently, the claim follows from \eqref{eq:neg} and (\ref{eq:double0}--\ref{eq:I_1}) by letting $S_2,S_1\rw\IF$.    \QED\\

Before proceeding with the proof of \netheo{th2} observe first that the variance function of $Y_u(s,t)$ is given by
\BQNY
\xx{V_Y}^2(s,t)=\frac{(1-\gamma)t^{2H}+(\Ga^2-\Ga)s^{2H}+\Ga(t-s)^{2H}}{(1+ct-c\Ga s)^2}.
\EQNY
In fact, the distribution function of $Y_u$ does not depend on $u$, so in the following we deal with $Y(s,t):=\frac{X_H(t)-\gamma X_H(s)}{1+ct-c\Ga s}$ instead of $Y_u(s,t)$. The \Hr{next} lemma \rE{\Hr{will be} used in} the proof of \netheo{th2}.

\BEL\label{lemVarY}
The variance function $\cJI{\xx{V^2_Y}}(s,t)$ attaints its unique global maximum over set $B:=\{(s,t): 0\le s\le t<\IF\}$ at $(\tilde{s}_0, \tilde{t}_0)$, with $\tilde{s}_0=0$
and $\tilde{t}_0=\frac{H}{c(1-H)}$. Further
$$
\xx{V_Y}(0, \tilde{t}_0)=\frac{H^H (1-H)^{1-H}}{c^H }.
$$
\EEL
\COM{
\prooflem{lemVarY} \cc{We first consider the case that $H\neq 1/2$.} It follows  that
\BQNY
\frac{\partial \xx{V_Y}^2}{\partial s}(s,t)=\frac{\partial \xx{V_Y}^2}{\partial t}(s,t)=0 &\Leftrightarrow &  s=t.
\EQNY
Thus the maximum of $\xx{V_Y}^2(s,t)$ over the set $B$ is \rE{attained} either on the line $s=t$ or on the line $s=0$.
\rE{In fact}, on the line $s=t$,
the maximum of $\xx{V_Y}^2(s,t)$ is \rE{attained} at $\left(\frac{H}{(1-\Ga)c(1-H)},\frac{H}{(1-\Ga)c(1-H)}\right)$, and on the line $s=0$,  the maximum of $\xx{V_Y}^2(s,t)$ is \rE{attained} at $\left(0,\frac{H}{c(1-H)}\right)$.
The claim follows by comparing the values on these two points.
\cc{Since $\cJI{\xx{V^2_Y}}(s,t)$ does  not have any extremal point in the case that $H=1/2$, the claim follows similarly, and thus the proof is complete.} \QED\\
}

\prooftheo{th2} The theorem will be proved in the following two steps.

\textbf{Step 1.} Let $K>\tilde{t}_0$ be a sufficiently large integer. We first \rE{derive} the asymptotics of
\BQNY
\pk{\sup_{0\le s \le t<K} Y(s,t)>u^{1-H}}\ \ \text{as}\ u\rw\IF.
\EQNY
Define  $\tilde{B}_{\tildel}:=\{(s,t): s\in(0,\tildel), t\in(\tilde{t}_0-\tildel,\tilde{t}_0+\tildel)\}$, for $\tildel>0$ sufficiently small, and let $B_K:=\{(s,t):0\le s\le t< K\}$. We write
\COM{%%%%%%%%%%%%%55555
\BQNY
\pk{\sup_{0\le s \le t<K} Y(s,t)>u^{1-H}}=\pk{\sup_{0\le s \le t<K} \overline{Y}(s,t)\frac{\xx{V_Y}(s,t)}{\xx{V_Y}(0,\tilde{t}_0)}>\frac{u^{1-H}}{\xx{V_Y}(0,\tilde{t}_0)}}
\EQNY
and}%%%%%%%%%%%%%555555
\BQNY
\pi(u):=\pk{\sup_{(s,t)\in \tilde{B}_{\tildel}}\tilde{Y}(s,t)>\frac{u^{1-H}}{\xx{V_Y}(0,\tilde{t}_0)}}\ \ \text{with} \ \tilde{Y}(s,t):=\overline{Y}(s,t)\frac{\xx{V_Y}(s,t)}{\xx{V_Y}(0,\tilde{t}_0)}.
\EQNY
Clearly
\BQNY
\pi(u)\le\pk{\sup_{0\le s \le t<K} Y(s,t)>u^{1-H}}\le \pi(u)+\pk{\sup_{(s,t)\in B_K/\tilde{B}_{\tildel}} \tilde{Y}(s,t)>\frac{u^{1-H}}{\xx{V_Y}(0,\tilde{t}_0)}}.
\EQNY
Therefore, we can conclude \yy{that}
\BQNY
\pk{\sup_{0\le s \le t<K} Y(s,t)>u^{1-H}}=\pi(u)\oo\ \ \text{as}\ u\rw\IF,
\EQNY
if
\BQN
\pk{\sup_{(s,t)\in B_K/\tilde{B}_{\tildel}} \tilde{Y}(s,t)>\frac{u^{1-H}}{\xx{V_Y}(0,\tilde{t}_0)}}=o(\pi(u))\ \ \text{as}\ u\rw\IF.\label{eqneg2}
\EQN
We \rE{calculate} next the aymptotics of $\pi(u)$ and show that \eqref{eqneg2} holds.
In view of \nelem{lemVarY} the standard \rE{deviation} function of $\tildeY$ given by
$$\sigma_{\tildeY}(s,t):=\frac{\xx{V_Y}(s,t)}{\xx{V_Y}(0,\tilde{t}_0)},\ \ (s,t)\in B,$$
 attains its unique maximum  over $\tilde{B}_{\tildel}$ at \yy{the point} $(0,\tildet)$, and $\sigma_{\tildeY}(0,\tildet)=1$.
\rE{Straightforward calculations yield }
\COM{\BQNY
(1+ct-c\gamma s)^{-1}=\frac{1}{1+c\tildet}\left(1+\frac{c}{1+c\tildet}(\tildet-t)+\frac{c\gamma}{1+c\tildet}s+\left(\frac{c}{1+c\tildet}\right)^2(\tildet-t+\gamma s)^2(1+o(1))\right)
\EQNY
and
\BQNY
\lefteqn{ ((1-\gamma)t^{2H}+(\gamma^2-\gamma)s^{2H}+\Ga(t-s)^{2H})^{1/2}}\\
&=&\tildet^{H}\Bigg(1-\frac{H}{\tildet}(\tildet-t)+\frac{\gamma^2-\gamma}{2\tildet^{2H}}s^{2H}-\frac{H\Ga}{\tildet}s
+\frac{(1-\Ga)H(2H-1)}{2\tildet^{2}}(\tildet-t)^2\\
&&
+\frac{\Ga H(2H-1)}{2\tildet^{2}}\left((\tildet-t)+s\right)^2 %+o(\left(\frac{1}{\tildet}(\tildet-t)-\frac{\gamma^2-2\gamma}{\tildet}s\right)^2)\\
-\frac{1}{8}\Big(\frac{2H}{\tildet}(\tildet-t)-\frac{\gamma^2-\gamma}{\tildet^{2H}}s^{2H}+\frac{2H\Ga}{\tildet}s
-\frac{(1-\Ga)H(2H-1)}{\tildet^{2}}(\tildet-t)^2\\
&&-\frac{\Ga H(2H-1)}{\tildet^{2}}\left((\tildet-t)+s\right)^2\Big)^2(1+o(1))\Bigg)
\EQNY
as $(s,t)\rw (0,\tildet)$.
%%%%%%%%%%555555555555555555555555555555555555555555555555555
\COM{
it follows that
\BQNY
1-\sigma_{\tildeY}(s,t)=\frac{c(\gamma-\gamma^2)}{2}s+\frac{c^2}{8}(\tildet-t)^2+\frac{c^2(3\gamma-2\gamma^2)}{4}(\tildet-t)s+o(.........)
\EQNY
Comparing this with \eqref{eqcov2}, we conclude that  Theorem 8.2 of Piterbarg (1996) can not be applied directly.
If $\gamma=0$, the explicit formula for the ruin probability is known, and it follows from Albrecher et al. (2008) we have
$$
\psi_{\Ga,\IF}(u)=1-(1-\psi_{0,\IF}(u))^{\frac{1}{1-\Ga}},
$$
for $H=1/2$. Unfortunately, even for this case we could not prove
$$
\psi_{\Ga,\IF}(u)={\frac{1}{1-\Ga}}\psi_{0,\IF}(u)(1+o(1)).
$$
It will be good if we can prove the last formula for all $H$.

-------------(the  following two formulas are wrong)-----------------

It can be seen that
\BQNY
\xx{V_Y}(s,t)=\left\{
              \begin{array}{ll}
\frac{\tilde{t}_0^H}{1+c\tilde{t}_0}\left(1-c^2(1-H)^2(\tildet-t)^2-c^2\Ga^2(1-H)^2 s^2\right)+o(\abs{\tildet-t}+s^2)    , &  H>1/2 ,\\
\frac{\tilde{t}_0^H}{1+c\tilde{t}_0}\left(1-c^2(1-H)^2(\tildet-t)^2-\frac{\Ga-\Ga^2}{2}{\tildet}^{-2H} s\right)+o(\abs{\tildet-t}+s), & H=1/2,\\
\frac{\tilde{t}_0^H}{1+c\tilde{t}_0}\left(1-c^2(1-H)^2(\tildet-t)^2-\frac{\Ga-\Ga^2}{2}{\tildet}^{-2H} s^{2H}\right)+o(\abs{\tildet-t}+s^{2H}), & H<1/2,
              \end{array}
            \right.
\EQNY
as $t\rw \tildet, s\rw 0$.
}
%%%%%%%%%%%%%%%%%%%%5555555555555555555555555555555555555555555555555555
Thus we conclude that
}
\BQN\label{eqv2}
1-\sigma_{\tildeY}(s,t)=\left\{
              \begin{array}{ll}
\frac{c^2(1-H)^3}{2H}(\tildet-t)^2(1+o(1))+\frac{(\Ga-\Ga^2)(1-H)^{2H}c^{2H}}{2H^{2H}} s^{2H}(1+o(1))    , &  H\le1/2 ,\\
%c^2(1-H)^2(\tildet-t)^2+\frac{\Ga-\Ga^2}{2}{\tildet}^{-2H} s+o(\abs{\tildet-t}+s), & H=1/2,\\
\frac{c^2(1-H)^3}{2H}(\tildet-t+\Ga s)^2(1+o(1))+\frac{(\Ga-\Ga^2)(1-H)^{2H}c^{2H}}{2H^{2H}} s^{2H}(1+o(1)), & H>1/2
              \end{array}
            \right.
\EQN
as $(s,t)\rw (0,\tildet)$. Additionally
\BQN\label{eqcov2}
1-Cov(\overline{Y}(s,t),\overline{Y}(s',t'))=\frac{1}{2\tildet^{2H}}\left(\mid t-t'\mid^{2H}+\Ga^{2}\mid s-s'\mid^{2H}\right)(1+o(1))
\EQN
as $(s,t), (s',t')\rw (0,\tildet)$ and for any $s,t,s',t'\in\tilde{B}_{\tildel}$
\BQNY
\E{\tildeY(s,t)-\tildeY(s',t')}^2
\COM{&=&\frac{1}{(\xx{V_Y}(0,\tilde{t}_0))^2}\E{Y(s,t)-Y(s',t')}^2\\
&\le&\frac{2}{(\xx{V_Y}(0,\tilde{t}_0))^2}\Bigg(\E{\frac{(X_H(t)-\Ga X_H(s))-(X(t')-\Ga X(s'))}{1+ct-c\Ga s}}^2\\
&+&\E{X(t')-\Ga X(s')}^2\left(\frac{1}{1+ct-c\Ga s}-\frac{1}{1+ct'-c\Ga s'}\right)^2\Bigg)\\
}
&\le& \CC (\abs{t-t'}^{2H}+\abs{s-s'}^{2H}).
\EQNY
Using \netheo{ThmPiter} for $H\le1/2$ and \netheo{ThmGPiter} \cc{with Remark \ref{b3n}} for $H>1/2$, we conclude that
\BQN
%\mathbb{P}\left(\sup_{t\in[0,T]}W_\Ga(t)>u\right)
\pi(u)=\mathcal{W}_{H}(u)\Psi\left(\frac{c^H u^{1-H}}{H^H (1-H)^{1-H}}\right)(1+o(1)),\label{Zvp1}
\EQN
where
$$
\mathcal{W}_{H}(u)=
 2^{\frac{1}{2}-\frac{1}{2H}}\frac{\sqrt{\pi}}{\sqrt{H(1-H)}} \mathcal{H}_{{2H}}\piter_{2H}^{\frac{1-\Ga}{\Ga}}\left(\frac{c^H u^{1-H}}{H^H (1-H)^{1-H}}\right)^{(1/H-1)}.
$$
Next we give the proof of \eqref{eqneg2}. Since $\sigma_{\tildeY}(s,t)$ is continuous, there exists some positive constant $\rho$ such that
  $$ \sup_{(s,t)\in B_K/\tilde{B}_{\tildel}} \sigma_{\tildeY}(s,t)<\rho<1  $$
  for the chosen small $\tildel$.
Therefore, in view of  Borell-TIS inequality (e.g., Adler and Taylor (2007)), for $u$ sufficiently large
\BQNY
\pk{\sup_{(s,t)\in B_K/\tilde{B}_{\tildel}} \tilde{Y}(s,t)>\frac{u^{1-H}}{\xx{V_Y}(0,\tilde{t}_0)}} \le \exp\left(-\frac{ (u^{1-H}-a)^2}{2 \xx{V_Y}^2(0,\tilde{t}_0)\rho^2}\right)
\EQNY
for some constant $a>0$. Consequently, Eq. \eqref{eqneg2} is established \cJ{ by comparing the last inequality with \eqref{Zvp1}}.\\

\textbf{Step 2.} We show that, for the chosen large enough integer $K>\tilde{t}_0$
\BQNY
\pk{\sup_{s,K \le t} Y(s,t)>u^{1-H}}=o\left(\pk{\sup_{0\le s \le t<K} Y(s,t)>u^{1-H}}\right)\asu.
\EQNY
\xH{For any $u>0$ we have (\cJI{set} $I_n=[n,n+1), n\inn$)} %It is derived that
\BQNY
&&\pk{\sup_{s,K \le t} Y(s,t)>u^{1-H}}
\le\pk{\sup_{s,K \le t} \frac{X_H(t)}{1+ct-c\Ga s}+\sup_{s,K \le t} \frac{-\Ga X_H(s)}{1+ct-c\Ga s}>u^{1-H}}\\
&&\le\pk{\sup_{s,K \le t} \frac{X_H(t)}{t^H}\frac{t^H}{1+ct-c\Ga s}>\frac{u^{1-H}}{2}}+\pk{\sup_{K\le s \le t} \frac{-\Ga X_H(s)}{s^H}\frac{s^H}{1+ct-c\Ga s}>\frac{u^{1-H}}{2}}\\
&&+\pk{\sup_{s\le K \le t} \frac{-\Ga X_H(s)}{(1+ct-c\Ga s)}>\frac{u^{1-H}}{2}}\\
%.
%\EQNY
%Let $I_n=[n,n+1), n\inn$. We further have %In view of Bonferroni's inequality
%\BQNY
%&&\pk{\sup_{s,K \le t} Y(s,t)>u^{1-H}}\\
%&&\le\pk{\sup_{s,K \le t} \frac{X_H(t)}{t^H}\frac{t^H}{1+ct-c\Ga s}>\frac{u^{1-H}}{2}}+\pk{\sup_{K\le s \le t} \frac{-\Ga X_H(s)}{s^H}\frac{s^H}{1+ct-c\Ga s}>\frac{u^{1-H}}{2}}\\
%&&+\pk{\sup_{s\le K }  X_H(s)>\frac{1+c(1-\Ga)K}{2\Ga}u^{1-H}}\\
%&&\le\sum_{i\ge K}\left(\pk{\sup_{t\in I_i} \frac{X_H(t)}{t^H}>\frac{1+c(1-\Ga)i}{2 i^H}u^{1-H}}+\pk{\sup_{s\in I_i} %\frac{X_H(s)}{s^H}>\frac{1+c(1-\Ga)i}{2\Ga i^H}u^{1-H}}\right)\\
%&&+\pk{\sup_{s\le K }  X_H(s)>\frac{1+c(1-\Ga)K}{2\Ga}u^{1-H}}\\
&&\le2 J_1(u)+J_2(u),
\EQNY
where %, for any $u>0$
\BQNY
J_1(u):=\sum_{i\ge K}\pk{\sup_{s\in I_i} \frac{X_H(s)}{s^H}>\frac{1+c(1-\Ga)i}{2 i^H}u^{1-H}},\ \ J_2(u):=\pk{\sup_{s\le K }  X_H(s)>\frac{1+c(1-\Ga)K}{2\Ga}u^{1-H}}.
\EQNY
 Furthermore, it follows that, for any $s,t\in I_i, i\ge K$
 \BQNY
\Hr {\E{ \Biggl( \frac{X_H(t)}{t^H}-\frac{X_H(s)}{s^H}\Biggr)^2}}=\frac{2s^Ht^H-2\E{X_H(t)X_H(s)}}{s^H t^H}\le\frac{\abs{t-s}^{2H}}{i^{2H}}\le \abs{t-s}^{2H}.
 \EQNY
Using Fernique's Lemma (e.g., Leadbetter et al. \cL{(1983)}) for some \xx{absolute} positive constants $C_1, C_2$
\BQNY
\pk{\sup_{t\in I_i} \frac{X_H(t)}{t^H}>\frac{1+c(1-\Ga)i}{2 i^H}u^{1-H}}\le C_1\expon{-C_2i^{2(1-H)}u^{2(1-H)}}
\EQNY
from which we conclude that, for $K$ sufficiently large
\BQNY
J_1(u)&\le& \sum_{i\ge K}  C_1\expon{-C_2i^{2(1-H)}u^{2(1-H)}}. %\le C_1\int_{K-1}^{\IF}\expon{-C_2x^{2(1-H)}u^{2(1-H)}}dx\\
%&\le&\frac{C_3}{2(1-H)}(K-1)^{1-2(1-H)}u^{-{2(1-H)}}\expon{-u^{2(1-H)}(K-1)^{2(1-H)}}
\EQNY
%for some positive constant $C_3$.
In the light of \eqref{eq:DS} \cL{of \netheo{HPDS}}
we see that
\BQNY
J_2(u)=\mathcal{D}_H\left(\frac{1+c(1-\Ga)K}{2\Ga K^H}u^{1-H}\right)^{\frac{(1-2H)_+}{H}}\Psi\left(\frac{1+c(1-\Ga)K}{2\Ga K^H}u^{1-H}\right)\oo \asu.
\EQNY
%where $\mathcal{F}_{H}$ is equal to $\mathcal{H}_{{2H}}2^{-\frac{1}{2H}}H^{-1}K^{1- 2H}$ if $H <1/2$, 2 if $H=1/2$, and 1 if $H>1/2.$
\COM{
$$
\mathcal{F}_{H}=\left\{
            \begin{array}{ll}
\mathcal{H}_{{2H}}\frac{2^{-\frac{1}{2H}}}{HK^{2H-1}} , & \hbox{if } \dE{H} <1/2 ,\\
2 & \hbox{if } \dE{H} =1/2,\\
1&  \hbox{if } \dE{H} >1/2,
              \end{array}
            \right.
$$}
Consequently, for sufficiently large $K$
\BQNY
\pk{\sup_{s,K \le t} Y(s,t)>u^{1-H}}&\le&2 J_1(u)+J_2(u)=o\left(\pk{\sup_{0\le s \le t<K} Y(s,t)>u^{1-H}}\right)
\EQNY
\rE{as $u\to \IF$, hence} the proof is complete.  \QED\\

\prooftheo{th1}
 Without loss of generality, we give only the proof of the case $\Ga\in(0,1).$
 Firstly, we give the asymptotic expansion of the standard deviation function $\xx{V_Z}(s,t)$ \cJI{at} \rE{the} point $(0,T)$.
 \rE{It follows that}
 \BQN\label{eq:Vzst}
\xx{V_Z}(s,t)
\COM{=\left((1-\gamma)t^{2H}+(\Ga^2-\Ga)s^{2H}+\Ga(t-s)^{2H}\right)^{1/2}\nonr\\
&&=T^H\left((1-\gamma)\left(1-\frac{T-t}{T}\right)^{2H}+(\Ga^2-\Ga)T^{-2H}s^{2H}+\Ga\left(1-\frac{T-t}{T}-\frac{s}{T}\right)^{2H}\right)^{1/2}\nonr\\
%&&=T^H\Big((1-\gamma)\left(1-2H\frac{T-t}{T}+o((T-t))\right)+(\Ga^2-\Ga)T^{-2H}s^{2H}\nonr\\
%&&\ \ +\Ga\left(1-2H\frac{T-t}{T}-2H\frac{s}{T}+o((T-t)+s)\right)\Big)^{1/2}\nonr\\
&&=T^H\left(1-2H\frac{T-t}{T}+(\Ga^2-\Ga)T^{-2H}s^{2H}-2H\Ga\frac{s}{T}+o((T-t)+s)\right)^{1/2}\nonr\\
&&
}=\left\{
              \begin{array}{ll}
T^H\left(1-HT^{-1}(T-t)-H\Ga T^{-1} s\right)+o((T-t)+s)    , &  H>1/2 ,\\
T^{1/2}\left(1-\cc{\frac{1}{2}}T^{-1}(T-t)-(\frac{1}{2}\Ga T^{-1}+\frac{\Ga-\Ga^{2}}{2}T^{-1}) s\right)+o((T-t)+s), & H=1/2,\\
T^H\left(1-HT^{-1}(T-t)-\frac{\Ga-\Ga^{2}}{2}T^{-2H} s^{2H}\right)+o((T-t)+s^{2H}), & H<1/2,
              \end{array}
            \right.
\EQN
 as $(s, t) \rw (0,T)$, hence % We see from \eqref{eq:Vzst} that
  there exists a positive constant $\delta>0$ such that
\BQN\label{eq:det}
  \abs{t-T-\gamma s}\le \mathbb{C} (\xx{V_Z}(0,T)-\xx{V_Z}(s,t))
\EQN
uniformly in $A_\delta:=\{(s,t): (s,t)\in[0,\delta]\times[T-\delta,T]\}$. %with $\delta$ chosen as above .
Next,
we study the asymptotics \Hr{of} the \rE{supremum of the} \Hr{Gaussian} \cL{random}  field \Hr{defined} on $A_\delta$.
\Hr{Set below}
\BQNY
\xx{\nu_u}(s,t)=\frac{u+ct-c\Ga s}{\xx{V_Z}(s,t)}\mq{and} \Pi(u)=\pk{\sup_{(s,t)\in A_\delta}\overline{Z}(s,t)\frac{\xx{\nu_u}(0,T)}{\xx{\nu_u}(s,t)}>\xx{\nu_u}(0,T)}.
\EQNY
For any $u>0$  %By the
\BQN
\Pi(u)\le\pk{\sup_{(s,t)\in A}(Z(s,t)-(ct -c\gamma s))>u}\le \Pi(u)+\pk{\sup_{(s,t)\in A/A_\delta}\overline{Z}(s,t)\frac{\xx{\nu_u}(0,T)}{\xx{\nu_u}(s,t)}>\xx{\nu_u}(0,T)}.\label{piZ}
\EQN
Since
\BQN\label{xixi}
\frac{\xx{\nu_u}(0,T)}{\xx{\nu_u}(s,t)}=1-\frac{\xx{V_Z}(0,T)-\xx{V_Z}(s,t)}{\xx{V_Z}(0,T)}-\frac{(c(t-T)-c\Ga s)\xx{V_Z}(s,t)}{(u+ct-c\Ga s)\xx{V_Z}(0,T)},
\EQN
we have, in view of \eqref{eq:det}, for any $\vp\in(0,1)$,  and sufficiently large $u$
\BQNY
1-\frac{\xx{V_Z}(0,T)-\xx{V_Z}(s,t)}{\xx{V_Z}(0,T)}\le \frac{\xx{\nu_u}(0,T)}{\xx{\nu_u}(s,t)}\le 1-(1-\vp)\frac{\xx{V_Z}(0,T)-\xx{V_Z}(s,t)}{\xx{V_Z}(0,T)}
\EQNY
uniformly in $(s,t)\in A_\delta$.
%Note that here in order to chose a suitable $\delta$ we have to use the expansion \eqref{eqcov1}.
Consequently
\BQN
\pk{\sup_{(s,t)\in A_\delta}Z_0(s,t)>\xx{\nu_u}(0,T)}\le \Pi(u)\le \pk{\sup_{(s,t)\in A_\delta}Z_\vp(s,t)>\xx{\nu_u}(0,T)},\label{pivp}
\EQN
where \xH{the random field} $\{Z_\vp(s,t), s,t\ge0\}$ is defined as
\BQNY
Z_\vp(s,t):=\overline{Z}(s,t)\left(1-(1-\vp)\frac{\xx{V_Z}(0,T)-\xx{V_Z}(s,t)}{\xx{V_Z}(0,T)}\right),\ \  \vp\in[0,1).
\EQNY
Direct calculations show that the standard \rE{deviation} function $\SZV(s,t):=\sqrt{\E{(Z_\vp(s,t))^2}}$ attains its unique maximum over $A_\delta$ at $(0,T)$ with $\SZV(0,T)=1$.
Thus, in the light of \eqref{eq:Vzst}, we have
\BQN\label{eqcov1}
\SZV(s,t)%&=&1-(1-\vp)\frac{\xx{V_Z}(0,T)-\xx{V_Z}(s,t)}{\xx{V_Z}(0,T)}\nonumber\\
&=&\left\{
              \begin{array}{ll}
1-(1-\vp)\left(HT^{-1}(T-t)+H\Ga T^{-1} s\right)(1+o(1))    , &  H>1/2 ,\\
1-(1-\vp)\left(\cc{\frac{1}{2}}T^{-1}(T-t)+(\frac{1}{2}\Ga T^{-1}+\frac{\Ga-\Ga^{2}}{2}T^{-1}) s\right)(1+o(1)), & H=1/2,\\
1-(1-\vp)\left(HT^{-1}(T-t)+\frac{\Ga-\Ga^{2}}{2}T^{-2H} s^{2H}\right)(1+o(1)), & H<1/2
              \end{array}
            \right.
\EQN
as $(s,t)\rw(0,T)$. Furthermore, it follows that
\BQN\label{cove1}
1-Cov(Z_\vp(s,t),Z_\vp(s',t'))%&=&1-\frac{\E{Z(s,t)Z(s',t')}}{\xx{V_Z}(s,t)\xx{V_Z}(s',t')}\nonumber\\
&=&\frac{1}{2T^{2H}}\left(\mid t-t'\mid^{2H}+\Ga^{2}\mid s-s'\mid^{2H}\right)(1+o(1))
\EQN
as $(s,t), (s',t')\rw(0,T)$. In addition, we obtain
\BQNY
\E{(Z_\vp(s,t)-Z_\vp(s',t'))^2}
%&=&\E{\vp\left(\frac{Z(s,t)}{\xx{V_Z}(s,t)}-\frac{Z(s',t')}{\xx{V_Z}(s',t')}\right)+(1-\vp)\frac{Z(s,t)-Z(s',t')}{\xx{V_Z}(0,T)}}^2\\
%&\le&2\vp^2\E{\frac{Z(s,t)}{\xx{V_Z}(s,t)}-\frac{Z(s',t')}{\xx{V_Z}(s',t')}}^2+\frac{2(1-\vp)^2}{\xx{V_Z}^2(0,T)}\E{Z(s,t)-Z(s',t')}^2\\
%&\le&\left(\frac{2\vp^2}{\xx{V_Z}^2(\delta,T-\delta)}+\frac{2(1-\vp)^2}{\xx{V_Z}^2(0,T)}\right) \E{Z(s,t)-Z(s',t')}^2\\
&\le& \CC (2|t-t'|^{2H}+2\Ga^2|s-s'|^{2H})
%&\le& 4C_{\delta,\vp} (1+\Ga^2)\left(\sqrt{(t-t')^2+(s-s')^2}\right)^{2H},
\EQNY
for $(s,t), (s',t')\in A_\delta$,  % and some positive constant $C_{\delta,\vp}$.
consequently, by \netheo{ThmPiter}
\BQN
\pk{\sup_{(s,t)\in A_\delta}Z_\vp(s,t)>\xx{\nu_u}(0,T)}=\mathcal{D}_{H,\Ga,\vp}\left(\frac{u+cT}{T^{H}}\right)^{(\frac{1}{H}-2)_+}\Psi\left(\frac{u+cT}{T^{H}}\right)(1+o(1))\label{Zvp}
\EQN
 as $u\rw\infty$, where
$$
\mathcal{D}_{H,\Ga,\vp}=\left\{
            \begin{array}{ll}
(1-\vp)^{-1}2^{-\frac{1}{2H}}H^{-1}\mathcal{H}_{{2H}}\piter_{2H}^{(1-\vp)^{\frac{1}{2H}}\frac{1-\Ga}{\Ga}}  , & \hbox{if } \dE{H} <1/2 ,\\
\piter_{1}^{(1-\vp)\frac{2-\Ga}{\Ga}}\times\piter_{1}^{(1-\vp)}[-\IF,0], & \hbox{if } \dE{H} =1/2,\\
1&  \hbox{if } \dE{H} >1/2
              \end{array}
            \right.
$$
and thus letting $\ve\rw0$, we obtain the asymptotic upper bound for $\Pi(u)$ on the set $A_{\delta}$. The asymptotic lower \cJI{bound} can be \Hr{derived} using the same arguments. In order to complete the proof we need to show further that
\BQN
\pk{\sup_{(s,t)\in A/A_\delta}\overline{Z}(s,t)\frac{\xx{\nu_u}(0,T)}{\xx{\nu_u}(s,t)}>\xx{\nu_u}(0,T)}=o(\Pi(u))\ \ \text{as}\ u\rw\IF.
\label{ZA}
\EQN
In the light of \eqref{xixi} for all $u$ sufficiently large
\BQNY
\sup_{(s,t)\in A/A_\delta} Var\left(\overline{Z}(s,t)\frac{\xx{\nu_u}(0,T)}{\xx{\nu_u}(s,t)}\right)\le (\rho(\delta))^2<1,
\EQNY
 \COM{ %%%%%%%%%%%%%%%%%%%%%%%%%%%%5555555555555555555555
 where $\rho(\delta)$ is a positive function in $\delta$ which exists due to the continuity of $\xx{V_Z}(s,t)$ in $A$.
Thus, in view of by Borell-TIS inequality, for $u$ sufficiently large
 Additionally, by the almost surely continuity of the process, for some constant $a>0$ and sufficiently large $u$ we have
\BQNY
\pk{\sup_{(s,t)\in A/A_\delta}\overline{Z}(s,t)\frac{\xx{\nu_u}(0,T)}{\xx{\nu_u}(s,t)}>a}\le\pk{\sup_{(s,t)\in A/A_\delta}\overline{Z}(s,t)\left(1-\frac{\xx{V_Z}(0,T)-\xx{V_Z}(s,t)}{2\xx{V_Z}(0,T)}\right)>a}   \le 1/2.
\EQNY
Therefore, direct application of the Borel inequality (e.g., Theorem D.1 of Piterbarg (1996)) implies
\BQNY
\pk{\sup_{(s,t)\in A/A_\delta}\overline{Z}(s,t)\frac{\xx{\nu_u}(0,T)}{\xx{\nu_u}(s,t)}>\xx{\nu_u}(0,T)}\le  \exp\left(\frac{(\xx{\nu_u}(0,T)-a)^2}{\rho^2(\delta)}\right)=o(\Pi(u)), \quad u\rw\IF.
\EQNY
for some constant $a>0$.
Consequently, Eq. (\ref{ZA}) is established, and thus the proof is complete.
}%55555555555555555555555555555555555555555555555555555555555555555555555555
where $\rho(\delta)$ is a positive function in $\delta$ which exists due to the continuity of $\xx{V_Z}(s,t)$ in $A$. Additionally, by the almost surely continuity of the \cJI{random field}, we have, for some constant $a>0$ %, and sufficiently large $u$
\BQNY
\pk{\sup_{(s,t)\in A/A_\delta}\overline{Z}(s,t)\frac{\xx{\nu_u}(0,T)}{\xx{\nu_u}(s,t)}>a}\le\pk{\sup_{(s,t)\in A/A_\delta}\overline{Z}(s,t)\left(1-\frac{\xx{V_Z}(0,T)-\xx{V_Z}(s,t)}{2\xx{V_Z}(0,T)}\right)>a}   \le 1/2.
\EQNY
Therefore, a direct application of the Borell inequality (e.g., Theorem D.1 of Piterbarg (1996)) implies
\BQNY
\pk{\sup_{(s,t)\in A/A_\delta}\overline{Z}(s,t)\frac{\xx{\nu_u}(0,T)}{\xx{\nu_u}(s,t)}>\xx{\nu_u}(0,T)}\le 2 \Psi\left(\frac{\xx{\nu_u}(0,T)-a}{\rho(\delta)}\right)=o(\Pi(u)) \ \ \text{as}\ u\rw\IF.
\EQNY
Consequently, Eq. (\ref{ZA}) is established, and thus the proof is complete. \QED\\

\section{Appendix}
The next theorem consists of two known results given in
H\"{u}sler and Piterbarg (1999) for the case $T=\IF$ and in D\c{e}bicki and Rolski (2002) when \cc{$T\in(0,\IF)$}.

\BT\label{HPDS}
\rE{If} $\{X_H(t),t\ge0\}$ is a fBm with Hurst \rE{index} $H\in(0,1]$, then for any $H\in(0,1)$
\BQN
\psi_{0,\IF}(u)=2^{\frac{1}{2}-\frac{1}{2H}}\frac{\sqrt{\pi}}{\sqrt{H(1-H)}} \mathcal{H}_{{2H}}\left(\frac{c^H u^{1-H}}{H^H (1-H)^{1-H}}\right)^{1/H-1}\Psi\left(\frac{c^H u^{1-H}}{H^H (1-H)^{1-H}}\right)(1+o(1))\label{eq:HP}
\EQN
holds as $ u\rw\infty$, and for any $H\in(0,1]$ \yy{and $T\in (0,\IF)$}
\BQN
\psi_{0,T}(u)=\mathcal{D}_H\left(\frac{u+cT}{T^H}\right)^{\frac{(1-2H)_+}{H}}\Psi\left(\frac{u+cT}{T^H}\right)(1+o(1))\label{eq:DS}
\EQN
holds as $ u\rw\infty$, where $\mathcal{D}_H$ is equal to $H^{-1}2^{-1/(2H)}\mathcal{H}_{2H}$ if $H<1/2$, 2 if $H=1/2$, and 1 if $H>1/2$.
\ET

 \cc{In the following theorem we present some \rE{results} used in the proof of our main theorems; denote the Euler Gamma function by $\Gamma(\cdot)$.}

\BT\label{ThmPiter}
Let $S,T$ be two positive constants, and \xH{let} $\{X(s,t),(s,t)\in\lbrack0,S]\times\lbrack0,T]\}$ be a zero-mean \cJ{Gaussian random field} with standard \cc{deviation} function $\sigma(\cdot,\cdot)$ and correlation function $r(\cdot,\cdot,\cdot,\cdot)$. Assume \xH{that} $\sigma(\cdot,\cdot)$ \cc{attains its unique maximum \xx{on} $\lbrack0,S]\times\lbrack0,T]$ at  $(s_0,t_0)$}, \cc{and further}
\BQN\label{PF}
\sigma(s,t)=1-b_{1}|s-s_{0}|^{\beta_{1}}(1+o(1))-b_{2}|t-t_{0}%
|^{\beta_{2}}(1+o(1)),\ \ \text{as }(s,t)\rightarrow(s_{0},t_{0})
\EQN
for some positive constants $b_i,\beta_i, i=1,2.$
Let, moreover %for correlation function be fulfilled
\[
r(s,s^{\prime},t,t^{\prime})=1-(a_1|s-s^{\prime}|^{\alpha_{1}}%
+a_2|t-t^{\prime}|^{\alpha_{2}})(1+o(1))\text{\ \ \ as }(s,t),(s^{\prime
},t^{\prime})\rightarrow(s_{0},t_{0})
\]
for some positive constants \cc{$a_i, i=1,2$ and $\alpha_i\in(0,2], i=1,2.$}
In addition, there \cL{exist two positive constants $G, \mu$ with $\mu\in(0,2]$} such that
$$
\E{(X(s,t)-X(s',t'))^2}\le G(|s-s^{\prime}|^{\xx{\mu}}%
+|t-t^{\prime}|^{\xx{\mu}})
$$
for any $(s,t), (s',t')\in\lbrack0,S]\times\lbrack0,T]$.
Then \rE{as $u\to \IF$}

i) if  $ \alpha_1<\beta_1$ and $ \alpha_2<\beta_2$
\[
\pk{\sup_{(s,t)\in\lbrack0,S]\times\lbrack0,T]}X(s,t)>u}= \prod_{i=1}^2 \Biggl( \mathcal{H}_{\alpha_i}a_i^{1/\alpha_i} b_i^{-1/{\beta_i}}\ \widehat{I}_i \  \Gamma\left(\frac{1}{\beta_i}+1\right) u^{2/\alpha_i-2/\beta_i}\Biggr)\Psi(u)(1+o(1));
\]

ii) if  $ \alpha_1<\beta_1$ and $ \alpha_2=\beta_2$
\[
\pk{\sup_{(s,t)\in\lbrack0,S]\times\lbrack0,T]}X(s,t)>u}= \mathcal{H}_{\alpha_1} a_1^{1/\alpha_1} b_1^{-1/{\beta_i}}\ \widehat{I}_1 \  \Gamma\left(\frac{1}{\beta_1}+1\right) \ \widehat{\mathcal{P}_{\alpha_2}^{b_2/ a_2}} u^{2/\alpha_1-2/\beta_1}\Psi(u)(1+o(1));
\]

iii) if  $ \alpha_1=\beta_1$ and $ \alpha_2=\beta_2$
\[
\pk{\sup_{(s,t)\in\lbrack0,S]\times\lbrack0,T]}X(s,t)>u}= \widehat{\mathcal{P}_{\alpha_1}^{b_1/a_1}}\ \widehat{\mathcal{P}_{\alpha_2}^{b_2/a_2}}\Psi(u)(1+o(1));
\]

iv) if  $ \alpha_1>\beta_1$ and $ \alpha_2>\beta_2$
\[
\pk{\sup_{(s,t)\in\lbrack0,S]\times\lbrack0,T]}X(s,t)>u}=\Psi(u)(1+o(1)),
\]
where \cc{ $\widehat{I}_2$ is the same as in \eqref{eq:main}  and
\BQNY
\widehat{\mathcal{P}_{\alpha_1}^{b_1/a_1}}:=\left\{
            \begin{array}{ll}
      \widetilde{\mathcal{P}}_{\alpha_1}^{b_1/a_1}      , & \hbox{if } s_0\in(0,S) ,\\
\mathcal{P}_{\alpha_1}^{b_1/a_1}, & \hbox{if } s_0=0\ \text{or}\ S,
              \end{array}
            \right.\ \
            \widehat{\mathcal{P}_{\alpha_2}^{b_2/a_2}}:=\left\{
            \begin{array}{ll}
      \widetilde{\mathcal{P}}_{\alpha_2}^{b_2/a_2}      , & \hbox{if } t_0\in(0,T) ,\\
\mathcal{P}_{\alpha_2}^{b_2/a_2}, & \hbox{if } t_0=0\ \text{or}\ T,
              \end{array}
            \right.\ \
              \widehat{I}_1:=\left\{
            \begin{array}{ll}
  2    , & \hbox{if } s_0\in(0,S) ,\\
1, & \hbox{if } s_0=0\ \text{or}\ S.
              \end{array}
            \right.
\EQNY
 }
\ET
\prooftheo{ThmPiter} In the context of Piterbarg (1996) and Fatalov (1992)  condition \eqref{PF} is formulated as
\BQN\label{PF2}
\sigma(s,t)=1-(b_{1}|s-s_{0}|^{\beta_{1}}+b_{2}|t-t_{0}%
|^{\beta_{2}})(1+o(1)),\ \ \text{as }(s,t)\rightarrow(s_{0},t_{0}).
\EQN
In fact, conditions \eqref{PF} and \eqref{PF2} play the same roles in the proof, since only bounds of the form
\BQNY
(b_{1}|s-s_{0}|^{\beta_{1}}+b_{2}|t-t_{0}%
|^{\beta_{2}})(1-\epsilon)\le 1-\sigma(s,t)\le (b_{1}|s-s_{0}|^{\beta_{1}}+b_{2}|t-t_{0}%
|^{\beta_{2}})(1+\epsilon)\
\EQNY
for any $\epsilon>0, \ \text{as }(s,t)\rightarrow(s_{0},t_{0}),$ are needed.
Therefore, the claims follow by  similar argumentations as in Piterbarg (1996) and Fatalov (1992). \QED\\

\COM{
\begin{remarks}

1, In the context of Piterbarg (1996) and Fatalov (1992), the condition \eqref{PF} is formulated as
\BQN\label{PF2}
\sigma(s,t)=1-(b_{1}|s-s_{0}|^{\beta_{1}}+b_{2}|t-t_{0}%
|^{\beta_{2}})(1+o(1)),\ \ \text{as }(s,t)\rightarrow(s_{0},t_{0}).
\EQN
In fact, the conditions \eqref{PF} and \eqref{PF2} play the same roles, since only bounds of the form
\BQNY
(b_{1}|s-s_{0}|^{\beta_{1}}+b_{2}|t-t_{0}%
|^{\beta_{2}})(1-\epsilon)\le 1-\sigma(s,t)\le (b_{1}|s-s_{0}|^{\beta_{1}}+b_{2}|t-t_{0}%
|^{\beta_{2}})(1+\epsilon),\
\EQNY
for any $\epsilon>0, \ \text{as }(s,t)\rightarrow(s_{0},t_{0}),$ are used in the proofs.\\
\end{remarks}}

\yy{{\bf Acknowledgement}: E. Hashorva and L. Ji kindly acknowledge partial
support from Swiss National Science Foundation Project 200021-1401633/1. The research of V.I. Piterbarg is supported by Russian Foundation for Basic Research, Project 11-01-00050-a. All the authors kindly acknowledge partial support by the project RARE -318984, a Marie Curie International Research Staff Exchange Scheme Fellowship within the 7th European Community Framework Programme}.

\end{document}